\documentclass[12pt]{article}
\usepackage{amssymb,latexsym}
\usepackage{fullpage}
\usepackage{array}
\renewcommand{\thefootnote}{}

\begin{document}
\title{A coincidence theorem for holomorphic maps to $G/P$}
\author{Parameswaran Sankaran    \\  
Institute of Mathematical Sciences\\
CIT Campus, Chennai 600 113, INDIA \\
E-mail: {\tt sankaran@imsc.ernet.in}\\ [4mm]} 
\date{}
\maketitle

\footnote{2000 A.M.S. Subject Classification:-  32H02, 54M20.}

\thispagestyle{empty}

\def\theequation {\arabic{section}.\arabic{equation}}
\renewcommand{\thefootnote}{}

\newcommand{\codim}{\mbox{{\rm codim}$\,$}}
\newcommand{\stab}{\mbox{{\rm stab}$\,$}}
\newcommand{\lr}{\mbox{$\longrightarrow$}}

\newcommand{\cm}{{\cal M}}
\newcommand{\cf}{{\cal F}}
\newcommand{\cd}{{\cal D}}

\newcommand{\blr}{\Big \longrightarrow}
\newcommand{\da}{\Big \downarrow}
\newcommand{\ua}{\Big \uparrow}
\newcommand{\hra}{\mbox{\LARGE{$\hookrightarrow$}}}
\newcommand{\rt}{\mbox{\Large{$\rightarrowtail$}}}
\newcommand{\dua}{\begin{array}[t]{c}
\Big\uparrow \\ [-4mm]
\scriptscriptstyle \wedge \end{array}}

\newcommand{\be}{\begin{equation}}
\newcommand{\ee}{\end{equation}}

\newtheorem{guess}{Theorem}
\newcommand{\bth}{\begin{guess}$\!\!\!${\bf .}~}
\newcommand{\eeth}{\end{guess}}
\renewcommand{\bar}{\overline}
\newtheorem{propo}[guess]{Proposition}
\newcommand{\bpropo}{\begin{propo}$\!\!\!${\bf .}~}
\newcommand{\epropo}{\end{propo}}

\newtheorem{lema}[guess]{Lemma}
\newcommand{\blem}{\begin{lema}$\!\!\!${\bf .}~}
\newcommand{\elem}{\end{lema}}

\newtheorem{defe}[guess]{Definition}
\newcommand{\bdefe}{\begin{defe}$\!\!\!${\bf .}~}
\newcommand{\edefe}{\end{defe}}

\newtheorem{coro}[guess]{Corollary}
\newcommand{\bcor}{\begin{coro}$\!\!\!${\bf .}~}
\newcommand{\ecor}{\end{coro}}

\newtheorem{rema}[guess]{Remark}
\newcommand{\brem}{\begin{rema}$\!\!\!${\bf .}~\rm}
\newcommand{\erem}{\end{rema}}

\newtheorem{exam}[guess]{Example}
\newcommand{\beg}{\begin{exam}$\!\!\!${\bf .}~\rm}
\newcommand{\eeg}{\end{exam}}

\newcommand{\ctext}[1]{\makebox(0,0){#1}}
\setlength{\unitlength}{0.1mm}
\newcommand{\cl}{{\cal L}}
\newcommand{\cp}{{\cal P}}
\newcommand{\cu}{{\cal U}}
\newcommand{\bz}{\mathbb{Z}}

\newcommand{\bq}{\mathbb{Q}}
\newcommand{\bt}{\mathbb{T}}
\newcommand{\bh}{\mathbb{H}}
\newcommand{\br}{\mathbb{R}}
\newcommand{\wt}{\widetilde}
\newcommand{\im}{{\rm Im}\,}
\newcommand{\bc}{\mathbb{C}}
\newcommand{\bp}{\mathbb{P}}
\newcommand{\spin}{{\rm Spin}\,}
\newcommand{\ds}{\displaystyle}
\newcommand{\tor}{{\rm Tor}\,}
\newcommand{\bs}{\mathbb{S}}
\def\ns{\mathop{\lr}}
\def\nssup{\mathop{\lr\,sup}}
\def\nsinf{\mathop{\lr\,inf}}
\renewcommand{\phi}{\varphi}
\newcommand{\co}{{\cal O}}
\begin{center}
{\it Dedicated to Professor Peter Zvengrowski}\\ 
{\it on the occasion of his sixty-first birthday}
\end{center}
\noindent
{\bf Abstract:} The purpose of this note is to extend 
to an arbitrary generalized Hopf and Calabi-Eckmann manifold the 
following result of Kalyan Mukherjea: Let $V_n=\bs^{2n+1}\times \bs^{2n+1}$ 
denote a Calabi-Eckmann manifold. If $f,g:V_n \lr \mathbb{P}^n$ are 
any two holomorphic maps, at least one of them being non-constant, then 
there exists a coincidence: $f(x)=g(x)$ for some $x\in V_n$. 
Our proof involves a coincidence theorem for holomorphic maps to    
complex projective varieties of the form $G/P$ where $G$ is complex 
simple algebraic group and $P\subset G$ is a maximal parabolic subgroup,
where one of the maps is dominant.
   
\section{Introduction}
Let $G$ be a simply connected simple algebraic group over $\bc$  
and  $P\subset G$ a {\it maximal} parabolic subgroup. 
Let $\cl$ be  the ample generator of $Pic(G/P)\cong \mathbb{Z}$. 
Let $E$ denote the total space of the principal $\bc^*$ bundle 
associated to $\cl^{-1}$.  Let $\lambda$ be any complex number with 
$|\lambda|>1$ and let $\phi:E\lr E$ denote the bundle map $e\mapsto 
\lambda\cdot e$, $e\in E$.  The quotient space, denoted $V_\lambda$, 
is a compact complex homogeneous non-K\"ahler manifold. The manifold 
$V_\lambda$ (or simply 
$V$) is called a generalized Hopf manifold \cite{rs}. (See also \S2,  
\cite{le}.)  One has an elliptic curve bundle  $q: V\lr G/P$ with  fibre and 
structure group the  elliptic curve $\bt=
\bc^*/\langle \lambda\rangle$ with periods $\{1, \tau\}$ where $\exp(2\pi 
\sqrt{-1}\tau)=\lambda$.  (Note that $\im(\tau)\neq 0$ as 
$|\lambda|>1$.) One has a diffeomorphism $V\cong \bs^1\times K/L$ where 
$K$ is a maximal compact subgroup of $G$ and $L$ 
the semi simple part of the centralizer of a subgroup of $K$ 
isomorphic to the circle $\bs^1$. If we take $G/P$ to be the complex 
projective space $\bp^n$, then the above construction yields 
the usual Hopf manifolds $\bs^1\times \bs^{2n+1}$.  

Let $G/P, G'/P'$, where both $G, G'$ are simply connected 
simple algebraic groups over $\bc$  
and $P, P'$, maximal parabolic subgroups of $G,G'$ respectively. Let
$E\lr G/P, E'\lr G'/P'$ denote 
the principal $\bc^*$-bundles associated to the negative ample 
generators of the Picard groups of $G/P,G'/P'$ respectively.  The
product bundle  $E\times E'\lr G/P\times G'/P'$ is a 
principal $\bc^*\times \bc^*$ bundle.  Let $\tau$ be any complex number 
with $Im(\tau)\neq 0$.  One has a complex analytic monomorphism  
$\bc\lr  \bc^*\times\bc^*$ 
defined by $z\mapsto (\exp(2\pi\sqrt{-1}\tau z),\exp(2\pi\sqrt{-1}z))$.  
Denote the image of this group by $\bc_\tau$.   
The group  $\bc_\tau\subset (\bc^*)^2$ acts on $E\times E'$ as bundle 
automorphisms and  
the quotient $U=E\times E'/\bc_\tau$ is the total space of an 
elliptic curve bundle with fibre and structure group the elliptic curve 
$(\bc^*)^2/\bc_\tau=\bt$ with periods $\{1,\tau\}$.  
Up to a diffeomorphism, 
$U$ can be identified with the space $K/L\times K'/L'$ where 
$K\subset G, K'\subset G'$ are maximal compact subgroups and $L$ and 
$L'$ are semi simple parts of centralizers of certain subgroups of $K$ and 
$K'$ isomorphic to $\bs^1$.  The compact complex manifold $U$ is homogeneous 
and non-K\"ahler, which we call a generalized Calabi-Eckmann manifold.  
When $G/P=\bp^m$, and $G'/P'=\bp^n$, $U$ is a Calabi-Eckmann 
manifold $\bs^{2m+1}\times \bs^{2n+1}$ \cite{ce}. We shall 
denote by $p$ the bundle projection $U\lr G/P\times G'/P'$. 

The manifold $U$ is an example of a simply connected compact complex 
homogeneous manifold. Such manifolds have been completely  
classified by H.-C. Wang \cite{w}. 

\bth \label{hopf} We keep the above notations.
(i) Let $\phi,\psi:V\lr G/P$ be any two holomorphic 
maps with $\phi$ non-constant. Then there exists an $x\in V$ such that 
$\phi(x)=\psi(x)$.  

(ii) Assume that $\dim(G/P)\leq \dim(G'/P')$ and 
let $\phi,\psi:U\lr G/P$ be a holomorphic map with $\phi$ 
non-constant. Then there exists an $x\in U$ such that $\phi(x)=\psi(x)$.  
\eeth

The above theorem will be derived from the following

\bth \label{singular}
(i) Let $M$ be any connected compact complex analytic manifold and let
$\phi, \psi:M\lr G/P$ be holomorphic where $P\subset G$ is a
maximal parabolic subgroup. Assume that at least one of the maps
$\phi,\psi$ is surjective.  Then there exists an $x\in M$
such that $\phi(x)=\psi(x)$.\\ 
(ii) Furthermore, if $M$ is a projective variety,  or if
K\"ahler manifold with $\dim(M)=\dim(G/P), $
and $f,g:M\lr G/P$ are continuous maps homotopic to $\phi, \psi$
respectively, then there exists an $x\in M$ such that $f(x)=g(x)$. 
\eeth

The special case of theorem \ref{hopf} when
$U=\bs^{2n+1}\times\bs^{2n+1}$ is a Calabi-Eckmann
manifold is due to K.Mukherjea \cite{m}. He 
uses D.Toledo's approach \cite{t} to analytic fixed point theory.  
In particular he uses Borel's computation of Dolbeault
cohomology of Calabi-Eckmann
manifolds and a certain `analytic Thom class' $\xi^0_\Delta\in 
H^n(\bp^n\times\bp^n;\Omega^n)\cong H^{2n}(\bp^n\times\bp^n;\bc)$ supported 
on the diagonal to detect coincidences.   
Our  proof is based on the observation that any holomorphic map 
from $U$  (resp. $V$) to a complex projective variety factors through  
$G/P\times G'/P'$ (resp. $G/P$) (see lemma \ref{mer}, \S3).  Any non-constant holomorphic map of $G/P$
or of $G/P\times G'/P'$ into $G/P$ is shown to be dominant 
(lemma \ref{picard}).
Theorem \ref{hopf} is then deduced from theorem \ref{singular}. 
Theorem \ref{singular} is proved using positivity of cup products in 
cohomology of $G/P$ and the fact that any effective cycle is
rationally equivalent to a {\it positive} linear combination of
Schubert cycles. (Cf. \cite{kl}, \cite{kuno}.)


\section{Holomorphic maps to $G/P$}
We keep the notations of \S 1. In particular, $G$ is a 
simply connected simple complex 
algebraic group.  Let  $Q\subset G$ be any parabolic subgroup
(not necessarily maximal).
Fix a maximal torus $T$ and a Borel subgroup
$B$ containing $T$ such that $B\subset Q$.   
Let $W$ be the Weyl group of $G$ with respect to $T$ and $W(Q)\subset W$ that 
of $Q$ and let $S\subset W$ denote the set of simple reflections 
with respect to our choice of $B$. Recall that the $T$ fixed points of 
$G/Q$ are labelled by $W/W(Q)$.
We shall identify $W/W(Q)$ with the set of coset representatives 
$W^Q\subset W$ having minimal length with respect to $S$. 
The $B$-orbits of these $T$ fixed points give an algebraic cell
decomposition for $G/Q$. 
Denote by $X(w)$ the Schubert variety
$X(w)\subset G/Q$ which is the $B$-orbit closure of the $T$-fixed point
corresponding to $w\in W^Q$. The class of the  
Schubert varieties  $[X(w)], w\in W^Q,$ form a $\bz$-basis for the
singular cohomology  group $H^*(G/Q;\bz)$ as well as the
Chow cohomology groups $A^*(G/Q)$. Recall that $Pic(G/Q)\cong A^1(G/Q)$
which is infinite cyclic when $Q$ is a maximal parabolic subgroup.   

\blem  \label{picard}
Let $P\subset G$ be a maximal parabolic subgroup and let $X$
be an irreducible complex projective variety with $Pic(X)=\bz$.
Suppose that $\dim(X)\geq \dim (G/P)$. Then any non-constant algebraic 
map $\phi: X\lr G/P$ is dominant with finite fibres; in particular 
$\dim(G/P)=\dim(X)$.  
\elem
\noindent
{\bf Proof:} Let $\cl$ be a very ample line bundle over $G/P$.
Let $Z=Im(\phi)$. As $\phi$ is non-constant, 
$Z$ is not a point variety and $\cl|Z$ is 
very ample. Since the bundle $\cl|Z$ is 
generated by global global sections, it follows that  
$\cl':=\phi^*(\cl)$ is generated by its sections over $X$. Also $\cl'$
cannot be trivial since any non-zero section of a trivial bundle 
is nowhere zero, whereas sections of $\cl'$ arising as pull-back 
of non-zero sections of $\cl|Z$ are non-zero sections which vanish 
somewhere in $X$. Since $Pic(X)\cong \bz$ is it follows that 
some {\it positive} multiple of $\cl'$ must be very ample. 
Hence, restricted to any fibres
of $\phi$ the bundle $\cl'=\phi^*(\cl)$ must be ample. This implies that  
the fibres of $\phi$ must be finite.  Therefore $\dim(X)
\leq \dim(G/P)$. Since $\dim(X)\geq \dim(G/P)$ by hypothesis,
we must actually have equality and the map $\phi$ must be dominant. 
\hfill $\Box$

\brem \label{picrem} (i) 
The assumption that $Pic(X)=\bz$ is not superfluous. For example, take 
$X=\bp^1\times \bp^n, ~n\geq 2$. Let $\phi:X\lr \bp^3$ be the composition
$\phi_1\circ pr_1,$ where $pr_1$ is the first projection map and  
$\phi_1:\bp^1 \lr \bp^3$ is defined by $\phi_1(z_0:z_1)=(z_0:z_1:0:0)$.\\
(ii)  Suppose that $\phi:Z\lr G/P$ is holomorphic 
and that $X\subset Z$ is a complex analytic subset  which satisfies the 
hypothesis of the lemma above. If $\phi:Z\lr G/P$ is holomorphic 
and $\phi|X$ is non-constant,  then $\phi$ must be dominant. On 
the other hand, suppose that $\pi:Z\lr M$ is a complex analytic 
fibre bundle with  fibre $X$ as in the lemma above. 
Suppose $\dim(X)>\dim(G/P)$, then 
any complex analytic map $\psi:Z\lr G/P$ factors through $\pi$, i.e, 
$\psi=\theta\circ\pi$ for some complex analytic map $\theta:M\lr G/P$
since the lemma implies that $\psi$ restricted to any fibre has to be  
constant.\\ 
(iii)  A result of K.Paranjape and V.Srinivas \cite{ps}
says that if  $G/P$ is not the projective space, any non-constant 
self morphism  $\phi:G/P \lr G/P$ is an automorphism of varieties.  
The full group of automorphisms of $G/P$ has been determined by 
I.Kantor \cite{k}. 
\erem

We shall now prove theorem \ref{singular}:\\

\noindent
{\bf Proof of theorem \ref{singular}:} (i)  
Suppose $\phi:M\lr G/P$ is surjective. Set $d=\dim(G/P), ~m=\dim(M)$. 
Let $\Delta\subset G/P\times G/P$ denote the diagonal of $G/P$.   
Let $\Gamma:=\Gamma_\theta\subset G/P\times G/P$ denote the image of  
the map $\theta:M\lr G/P \times G/P$, $\theta(x)=
(\phi(x),\psi(x)), ~x\in M$. 
We need only show that $\Gamma\cap\Delta\neq \emptyset.$ 
Note that $\Gamma$ is a complex analytic subspace of the 
projective variety $G/P$ 
and hence algebraic by GAGA \cite{serre}. 
Also $k:=\dim(\Gamma)\geq \dim(G/P)$ since $\phi$ is
surjective.  Now, as for any effective cycle in $G/P\times G/P$, 
the class $[\Gamma]\in A_k(G/P\times G/P)$ is a  
{\it positive} linear combination of Schubert cycles in $G/P\times G/P$. 
(Cf. \cite{k}.  See also  \cite{kuno}.) 
Thus, $[\Gamma]=\sum_i a_i [X(w_i)]\times [X(w'_i)]$ where $a_i$ 
are positive integers, $w_i, w'_i\in W^P$ are suitable elements such 
that $\dim(X(w_i))+\dim(X(w'_i))=\dim(\Gamma)$. Since $\phi$ is surjective, 
in the above expression for $\Gamma$, the 
term $[G/P]\times [X(w)]$ must occur with positive  
coefficient for some $w\in W^P$. The same arguments
can be applied to the class of the
diagonal $\Delta$ in $G/P\times G/P$ as well. 
In fact, it is known that $[\Delta]=
\sum_{v\in W^P} [X(v)]\times [X(v'')]$ where $X(v'')$ 
is the Schubert variety `dual' to $X(v)$, i.e, $v''=w_0.v$ where 
$w_0\in W^P$ represents the longest element of $W/W_P$.
(Cf. theorem 11.11, \cite{ms}.)
Hence  $[\Gamma].[\Delta]= (a[G/P]\times [X(w)]+other~terms)
\cdot(1\times [G/P]+other~terms)=a(1\times [X(w)])+other~terms$,
where $a>0$ and the coefficients of the remaining terms 
(with respect to the basis consisting of Schubert cocycles) in the
rhs of the last equality are {\it non-negative
integers.} Hence $[\Gamma][\Delta]\neq 0$ in $A_{k-d}(G/P\times G/P)$ 
and so $\Gamma\cap \Delta\neq\emptyset$.

\noindent 
(ii)  Let $h:M\lr G/P\times G/P$ be the map $x\mapsto (f(x),g(x))$ 
for $x\in M$.
It suffices to show that $h^*([\Delta])\in H^{2d}(M;\bz)$ is non-zero, where
we regard $[\Delta]$ as an element of the singular cohomology group 
$H^{2d}(G/P\times G/P;\bz)$. Note that $h$ is homotopic to $\theta
=(\phi,\psi)$ and so 
we have $h^*([\Delta])=\theta^*([\Delta])$. To complete the proof, 
it suffices to show that $\theta^*([\Delta])\neq 0$ in $H^*(M;\bz)$.

\noindent
Suppose that $M$ is K\"ahler and $\dim(M)=d=\dim(G/P)$. 
By de Rham and Hodge theory (\S 15.7, \cite{hirz}),  we have
$H^{r}(M;\bc)=\oplus_{p+q=r}H^{p,q}_{\bar{\partial}}(M)$ 
and $\theta^*([\Delta])$ can be 
thought of as an element of $H^{d,d}_{\bar{\partial}}(M;\bc)$. 
Since $\phi$ is dominant, $\dim(\Gamma)=d=\dim(M)$ 
and the fundamental class $\mu_M\in H_{2d}(M;\bz)$ maps to
$n[\Gamma]\in H_{2d}(G/P\times G/P;\bz)\cong A_d(G/P\times G/P)$
for some $n\geq 1$. From what has been shown already the intersection
product $[\Delta].[\Gamma_\theta]\neq 0$ in the Chow ring of $G/P\times G/P$.
This  implies that 
$[\Delta]\cap [\Gamma_\theta]\neq 0$ in $H_0(G/P\times G/P)$.
Now we have (see ch. 5, \S6, \cite{span}) \\
$$\theta_*(\theta^*([\Delta])\cap \mu_M)=n[\Delta]\cap \theta_*(\mu_M)
=n[\Delta]\cap [\Gamma]\neq 0.$$ 
It follows that $\theta^*([\Delta])\neq 0$ in $H^*(M;\bz)$.

\noindent
If $M$ is a complex projective variety then, one can always find 
an irreducible subvariety $Z\subset M$ with $\dim(Z)=\dim(\Gamma_\theta) $
which maps onto $\Gamma$. It follows that $\theta_*([Z])=
n[\Gamma_\theta]\in H_*(G/P\times G/P)$ for some $n>0$.  Using the 
fact that $H_*(G/P\times G/P)$ has no torsion, proceeding just as 
before we conclude that $\theta^*(\Delta)\neq 0$. 
\hfill $\Box$

\noindent
\brem (i) 
In the statement of theorem \ref{singular}, the 
hypothesis that $M$ be nonsingular is not necessary. 
Indeed, H.Hironaka \cite{hir} has shown that any irreducible 
complex analytic space which is countable at 
infinity  can be desingularized. So replacing $M$ by $\wt{M}$ and 
the maps $\phi,\psi$ by $\wt{\phi}:=\phi\circ \pi, \wt{\psi}:=\psi\circ \pi $ 
respectively where $\pi:\wt{M}\lr M$ is a desingularization map, the see that 
$\wt{\phi}, \wt{\psi}$ must have a coincidence. This immediately implies 
that $\phi$ and $\psi$ must have a coincidence. \\
\noindent
(ii) In case $M$ is not K\"ahler, it is not true in general that
$\theta^*([\Delta])\neq 0$ in $H^*(M)$ although $\phi$ and $\psi$ 
must have a coincidence as our theorem shows.  
Such coincidences 
have been described as ``homologically invisible" by Kalyan Mukherjee 
\cite{m}.  He observed that when $M$ is the Calabi-Eckmann manifold
$\bs^{2n+1}\times \bs^{2n+1}$ with $n\geq 1$, for any two maps $\phi ,\psi:M\lr \bp^n$ 
the homomorphism $\theta^*:H^*(\bp^n\times\bp^n)\lr H^*(M)$ is zero in 
positive dimensions where $\theta=(\phi,\psi)$. 
In particular, if $f, g$ are {\it continuous}  maps 
homotopic to holomorphic maps $\phi, \psi$ respectively with 
$\phi$ dominant, we do not know if  $f$ and $g$ must have a 
coincidence.\\
\noindent 
(iii) When $M$ is K\"ahler and  
$\dim(M)>\dim(G/P)$, the conclusion of the theorem is still valid provided 
$[\Gamma_\theta]\in H_*(G/P;\bq)$ is in the image of 
$\theta_*:H_*(M;\bq)\lr H_*(G/P;\bq)$.  
\erem

\noindent
\bcor \label{coingmodp}
(i)  Let  $P\subset G$ be a maximal parabolic subgroup and let
$f,g:G/P\lr G/P$ be any two continuous maps homotopic to
holomorphic maps $\phi, \psi$ respectively where $\phi$ is non-constant.
Then $f(x)=g(x)$ for some $x\in G/P$.

\noindent
(ii)  Let $f,g:G/P\times G'/P'\lr G/P$ be any two continuous maps which are
homotopic to holomorphic maps $\phi,\psi$ respectively with $\phi$  
being non-constant. Assume that $\dim(G/P)\leq \dim(G'/P')$.
Then there exists an $x\in G/P\times G'/P'$ such that 
$f(x)=g(x)$.
\ecor
\noindent
{\bf Proof:} Part (i) follows immediately from
lemma \ref{picard} and theorem \ref{singular}(ii).
To prove (ii), 
suppose  $\phi|Z$ is constant for every fibre 
$Z\cong G'/P'$ of the first projection
$pr_1:G/P\times G'/P'\lr G/P$ map, then $\phi$ can be factored as
$\phi_1\circ pr_1$ where $\phi_1:G/P\times G'/P'\lr G/P$ defined
by $\phi$.  It follows from \ref{picard} that $\phi_1$ 
is dominant. Hence $\phi$ is also dominant. 
Otherwise for some fibre $Z\cong G'/P'$, $\phi|Z$ is non-constant. By lemma 
\ref{picard} it follows that $\phi|Z$ is dominant. It follows from 
lemma \ref{picard} again 
that $\phi|Z$ --- and hence $\phi$ --- must be dominant and the corollary 
follows from theorem \ref{singular}. \hfill $\Box$

\brem\label{fpf}
It follows from the above corollary that any continuous map homotopic
to a holomorphic map has a fixed point. However,  
in general, the spaces $G/P$ do not have fixed point property. 
For example, the Grassmannian $G_k(\bc^n)
=SL(n,\bc)/P_k$ admits a continuous fixed point free involution 
whenever $n$ is even and $k$ odd or if $n=2k$. As another example, 
the complex quadric $SO(n)/P_1$ is diffeomorphic to the oriented 
real Grassmann manifold $\wt{G}_2(\br^n)$ of oriented $2$-planes in $\br^n$.
The involution that reverses the orientation on each element of 
$\wt{G}_2(\br^n)$ is obviously fixed point free.   
However, it is known that $G_k(\bc^n)$ has fixed point property 
(for continuous maps) when $n$ is large compared to $k$ and 
at most one of $n-k,k$ is odd.   See \cite{gh}, \cite{gh2}.  
\erem
 
\section {\bf Proof of Main Theorem} 

We now prove the main result of the paper, namely, theorem \ref{hopf}.
We keep the notations of \S 1. 

\blem\label{mer}
Let $U,V$ be generalized Calabi-Eckmann and generalized Hopf manifolds. 
(See \S1.)
Let $Z$ be a complex projective variety. 
Any holomorphic maps $\phi:U\lr Z$, $\psi:V\lr Z$ can be factored as 
$\phi=\phi_1\circ p, \psi=\psi_1\circ q$, where $p:U\lr G/P\times G'/P'$ and 
$q:V\lr G/P$ are projections of the principal $\bt$-bundles.
\elem

\noindent
{\it Proof:}  It was shown in the proof of Theorem 3, \cite{rs}, that  
$H^2(V;\bz)=H^2(K/L;\bz)=0$ where $K\subset G$ is a maximal compact subgroup 
of $G$ and $L$ is the semi simple part of the centralizer in $K$  of 
a subgroup of $K$ isomorphic to $\bs^1$.  
The same argument shows that $H^2(U;\bz)=
H^2(K/L\times K'/L';\bz)=0$.  In particular the manifolds $U,V$ are not 
K\"ahler.

A theorem of Grauert and Remmert \cite{gr} says that 
for a compact complex homogeneous manifold $M$ of dimension $n$, the 
transcendence degree over 
$\bc$ of the field ${\cal M}(M)$ of meromorphic functions on $M$ is 
equal to $n$ if and only if it a projective algebraic variety.
In the our case $U,V$ fibre over projective varieties of dimension 
$1$ less. It follows that $tr.deg_\bc({\cal}(V)) 
\geq tr.deg_\bc({\cal M}(X))=\dim(G/P)$.  
Suppose $\psi$ is not constant along a fibre.  Then there exists an 
open set (in the analytic topology)  $N\subset G/P$ such that 
$\psi$ is non-constant on $q^{-1}(x)$, for any  $x\in N$. Let $x\in N$.  
Composing with a suitable meromorphic function on $Z$ which is 
non-constant on $\psi(q^{-1}(x))$,  we get 
a meromorphic function $\theta$ on $V$. We  claim that 
$\theta$ is transcendental over ${\cal M}(G/P)\subset {\cal M}(V)$. Assume, 
if possible, that  
$\theta^k+a_1\theta^{k-1}+\cdots + a_k=0, a_i\in {\cal M}(G/P)$. By changing 
the $x\in N$ if necessary, we may assume that  $x$ is not 
on the polar divisor for any $a_i$. Restricting this equation to 
the fibre over $x$, we see that $\theta|q^{-1}(x)$ is algebraic over $\bc$. 
Since $\bc$ is algebraically closed, we must have $\theta|p^{-1}(x)\in \bc$.
This is absurd since $\theta$ is non-constant on $q^{-1}(x)$. 
Hence we conclude that $\theta$ is constant along the fibres of 
$q$. Proof that $\phi$ is constant along the fibres of 
$p$ is entirely similar.\hfill $\Box$

\noindent
{\bf Proof of Theorem \ref{hopf}:} (i) By lemma \ref{mer}, 
the maps $\phi, \psi$ factor through the 
projection of the elliptic curve bundle $p:V\lr G/P$. Write 
$\phi=\phi_1\circ p,~ \psi=\psi_1\circ p$. Now, it suffices to 
show that the holomorphic maps $\phi_1$ and $\psi_1$ have a coincidence.  
Since $\phi_1$ is non-constant, this is now immediate from 
corollary \ref{coingmodp} (i). \\
\noindent 
(ii) Proceeding exactly as in (i), we write $\phi=\phi_1\circ q, \psi
=\psi_2\circ q$, 
where $\phi_1,\psi_1:G/P\times G'/P'\lr G/P$ are holomorphic. Note that 
since $\phi_1$ is non-constant.  By corollary \ref{coingmodp} 
$\phi_1$ and $\psi_1$ must have a coincidence.  Hence $\phi$ and 
$\psi$ must also have a coincidence. \hfill $\Box$

We conclude with the following observation.

\noindent 
\blem\label{mer2}
Let $\pi:W\lr M$ be a holomorphic fibre bundle with $M$ 
compact connected, and fibre a complex torus $\bt$.
Suppose that $H_2(F;\bq)\lr H_2(W;\bq)$ is zero.  Then 
any holomorphic map $\phi:W\lr G/Q$ is constant on the 
fibres of $\pi$ where $Q\subset G$ is any parabolic subgroup.
\elem
\noindent
{\bf Proof:} Assume that $\phi:W\lr G/Q$ is a holomorphic 
map such that $\phi|F$ is not constant for some fibre 
$F$ of the $T$-bundle $\pi:W\lr M$. 
Let $\iota:F\subset W$ denote the inclusion map. 
Let $C\subset G/Q$ be the 
image of $F$.  Note that $\dim(C)=1=\dim (F)$.  Since 
$\phi$ is holomorphic, $C$ is an algebraic subvariety 
of $G/Q$. In particular, it represents a non-zero 
element of $H_2(G/Q;\bq)$. In fact $C$ is rationally equivalent to
a positive linear combination of certain $1$-dimensional Schubert
subvarieties in $G/Q$ \cite{kl}. It follows that $(\phi|F)_*:H_2(F;\bq)\lr
H_2(G/Q;\bq)$ maps the fundamental class of $F$ to a nonzero 
element of $H_2(G/Q;\bq)$. On the other hand, 
$(\phi|F)_*=\phi_*\circ \iota_*=0$ in dimension $2$, 
since $\iota_*:H_2(T;\bq)\lr H_2(W;\bq)$ is zero by hypothesis.
We conclude that $\phi$ must be constant on the fibres of $\pi$.
\hfill $\Box$

\noindent 
\brem
(i) Let $\pi:W\lr M$ be as in the above lemma. Let 
$\phi, \psi:W\lr G/P$ be any two holomorphic maps where $\phi$ is 
non-constant and $P\subset G$ a maximal parabolic. Let $\phi_1:M\lr G/P$ 
be such that $\phi=\phi_1\circ \pi$. 
Suppose $X\subset M$ is irreducible and has the structure of a   
complex projective variety with $Pic(X)=\bz$ and $\dim(X)=\dim(G/P)$. 
If $\phi_1|X$ is non-constant, then, in view of the above lemma and 
remark \ref{picrem}(ii), $\phi_1 $ must be dominant. It follows that 
$\phi$ itself must be dominant.  By theorem 
\ref{singular} it follows that  $\phi$ and $\psi$ 
must have a coincidence.  

(ii) I do not know if theorem \ref{hopf}  still holds if one 
merely assumes 
that $\phi,\psi$ are continuous maps homotopic to holomorphic maps one 
of which is dominant.   
\erem 

{\bf Acknowledgments:} I am indebted to  Kalyan Mukherjea for asking me 
the question which led to this paper and for pointing out his 
paper \cite{m}. 
I am grateful to the referee of this paper for his/her valuable 
comments and suggestions for improvements. 
I am grateful to Issai Kantor for a copy of his
paper \cite{k} and for translating into English the relevant parts 
of his of paper for my benefit.

\end{document}